\documentclass[12pt,oneside]{article}
\usepackage{amsmath}
\usepackage{amsfonts}
\usepackage{amsthm}
\usepackage{amssymb}
\usepackage{graphicx}
\usepackage{color}
\usepackage{txfonts}
\usepackage{setspace}
\usepackage{mathrsfs}
\usepackage[square, comma, sort&compress, numbers]{natbib}

\usepackage[colorlinks,
            citecolor=blue,
            linkcolor=blue]{hyperref}
\usepackage{xcolor}

\newtheorem{proposition}{Proposition}[section]

\newtheorem{definition}[proposition]{Definition}
\newtheorem{remark}[proposition]{Remark}
\newtheorem{theorem}[proposition]{Theorem}

\newtheorem{lemma}[proposition]{Lemma}
\numberwithin{equation}{section}

\begin{document}
\title{Stochastic Optimization Theory of Backward Stochastic Differential Equations Driven by G-Brownian Motion
\thanks{This work was Supported by National Basic Research Program of China [Project No.
973-2007CB814901].}
\author{Zhonghao Zheng$^1$ \ \ \ Xiuchun Bi$^2$\ \ \
  Shuguang Zhang$^2$\footnote{Corresponding author: sgzhang@ustc.edu.cn}
\\School of Mathematical Sciences, Xiamen University, \\Xiamen 361005 P.R. China
\\Department of Statistics and Finance, University of Science \\and Technology of
China, Hefei 230026  P.R. China}}
\date{}
%\fontsize{13}{0}\selectfont
\maketitle
\textbf{Abstract}:
In this paper, we consider the stochastic optimal control problems under G-expectation. Based on the theory of backward stochastic differential equations driven by G-Brownian motion, which was introduced in \cite{r20,r21}, we can investigate the more general stochastic optimal control problems under G-expectation than that were constructed in \cite{r22}. Then we obtain a generalized dynamic programming principle and  the value function is proved to be a viscosity solution of a fully nonlinear second-order partial differential equation.
\\\textbf{Keywords}:\ \ G-expectation, G-Brownian motion, Backward SDEs, Stochastic control.
\section{Introduction}
Non-linear BSDEs in the framework of linear expectation were introduced by Pardoux and Peng \cite{r24} in 1990. Then a lot of researches were studied by many authors and they  provided various applications of BSDEs in stochastic control, finance, stochastic differential games and second order partial differential equations theory, see \cite{r25,r26,r27,r28,r29,r30,r31,r32}.
\par
The notion of sublinear expectation space was introduced by Peng \cite{r1,r2,r3}, which is a generalization of classical probability space. The G-expectation, a type of sublinear expectation, has played an important role in the researches of sublinear expectation space recently. It can be regarded as a counterpart of the Wiener probability space in the linear case. Within this G-expectation framework, the G-Brownian motion is the canonical process. Besides, the notions of the G-martingales and the It\^{o} integral w.r.t. G-Brownian motion were also derived. There are some new structures in these notions and some new applications in the financial models with volatility uncertainty, see Peng \cite{r3,r4}.
\par
In the G-expectation framework, thanks to a series of studies \cite{r8,r9,r16,r17},  the complete representation theorem for G-martingales has been obtained by Peng, Song and Zhang \cite{r33}. Due to this contribution, a natural formulation of BSDEs driven by G-Brownian motion was found by Hu, Ji, Peng and Song \cite{r20}. In addition, the existence and uniqueness of the solution to the BSDEs driven by G-Brownian motion has been proved.  They also have given the comparison theorem, Feynman-Kac Formula and Girsanov transformation for  BSDEs driven by G-Brownian motion in \cite{r21}. So the complete theory of BSDEs driven by G-Brownian motion has been established.
\par
 An important application of BSDEs is that we can define the recursive utility functions from BSDEs, which can index scaling risks in the study of economics and finance \cite{r34,r35,r36,r37}. Based on these results, a type of significant stochastic optimal control problems under linear expectation with a BSDE as cost function were studied \cite{r28,r29,r30,r31,r32}. Under G-expectation, the similar problems will be useful in the future studies of  finance models with volatility uncertainty. So we arise a natural question: Can we construct the similar results in G-expectation framework. When the complete results about BSDEs driven by G-Brownian motion were established in \cite{r20,r21}, we try to prove the complete results of  stochastic optimization theory of BSDEs  driven by G-Brownian motion in this paper.
\par
In this paper, we investigate the stochastic optimal control problems with a BSDE driven by G-Brownian motion constructed in \cite{r20,r21} as cost function. Based on the results in \cite{r20,r21}, we obtain the dynamic programming principle under G-expectation. Besides, the value function is proved to be a viscosity solution of a fully nonlinear second-order partial differential equation.
\par
The rest of the paper is organized as follows. In Section 2, we recall the G-expectation framework and adapt it according to our objective. Besides, we give the related properties of forward and backward stochastic differential equations driven by G-Brownian motion, which will be needed in the sequel sections. In Section 3, the stochastic optimal control problems with a BSDE  driven by G-Brownian motion as cost function are investigated and a dynamic programming principle under G-expectation is obtained. In Section 4, The value function is proved to be a viscosity solution of a fully nonlinear second-order partial differential equation.

\section{Preliminaries}
In this section, we recall the  G-expectation framework established by Peng \cite{r1,r2,r3,r5}. Besides, we give some results about forward and backward stochastic differential equations driven by G-Brownian motion, which we need in the following sections. Some details can be found in \cite{r20,r21}.

\subsection{G-expectation and G-martingales}

 \begin{definition}
 Let $\Omega$ be a given set and $\mathcal{H}$ be a linear space of real valued functions defined on $\Omega$, namely  $c \in \mathcal{H}$ for each constant c and $\mid X\mid \in \mathcal{H}$ if $X \in \mathcal{H}$. The space $\mathcal{H}$ can be considered as the space of random variables. A sublinear expectation $\mathbb{E}$ is a functional $\mathbb{E}:\mathcal{H}\rightarrow \mathbb{R}$ satisfying the following properties: for all $X,Y \in \mathcal{H}$, we have

$\left(\rm{i}\right)$ Monotonicity: $\mathbb{E}[X]\geq \mathbb{E}[Y]$,\ if $X \geq Y$;

$\left(\rm{ii}\right)$ Constant preservation: $\mathbb{E}[c]=c$,\ for $c\in \mathbb{R}$;

$\left(\rm{iii}\right)$ Sub-additivity: $\mathbb{E}[X+Y]\leq \mathbb{E}[X]+\mathbb{E}[Y]$;

$\left(\rm{iv}\right)$ Positive homogeneity: $\mathbb{E}[\lambda X]=\lambda \mathbb{E}[X]$,  for $\lambda \geq 0$.
\par
The triple $(\Omega,\mathcal{H},\mathbb{E})$ is called a sublinear expectation space.
\end{definition}

\begin{definition}(G-Normal Distribution)  A d-dimensional random vector $X=(X_{1},\cdots,X_{d})$ on a sublinear expectation space $(\Omega, \mathcal{H}, \mathbb{E})$ is called G-normally distributed if for each $a,b\geq 0$, we have

\begin{equation*}
X+b\bar{X}\overset{d}{=}\sqrt{a^{2}+b^{2}}X,
\end{equation*}
where $\bar {X}$ is an independent copy of $X$, i.e., $\bar {X}$ and $X$ is identically distributed and $\bar{X}$ is independent from $X$. Here the letter $G$ denotes the function
\begin{equation}\label{10}
G(A):=\frac{1}{2}\mathbb{E}[(AX,X)]:\mathbb{S}_{d}\mapsto \mathbb{R},
\end{equation}
where $\mathbb{S}_{d}$ denotes the collection of all $d\times d$ symmetric matrices.
\end{definition}
\begin{proposition}
Let $X$ be G-normal distributed. The distribution of X is characterized by
\begin{equation}
u(t,x)=\mathbb{E}[\varphi(x+\sqrt{t}X)],\ \varphi \in C_{b.Lip}(\mathbb{R}^{d}).
\end{equation}
In particular, $\mathbb{E}[\varphi(X)]=u(1,0)$, where $u$ is the unique viscosity solution of the following parabolic PDE defined on $[0,\infty)\times \mathbb{R}^{d}$:
\begin{equation}
\partial_{t}u-G(D^{2}u)=0,\  u|_{t=0}=\varphi,
\end{equation}
where $G$ is defined by $\eqref{10}$.
\end{proposition}
\begin{remark}
It is easy to check that $G$ is a monotonic sublinear function defined on $\mathbb{S}(d)$ and $G(A):=\frac{1}{2}\mathbb{E}[(AX,X)]\leq \frac{1}{2}|A|\mathbb{E}[|X|^2]=\frac{1}{2}|A|\overline{\sigma}^{2}$ implies that there exists a bounded, convex and closed subset $\Gamma \subset \mathbb{S}^{+}_{d}$ such that
\begin{equation}\label{11}
G(A)=\frac{1}{2}\underset{\gamma \in \Gamma}{sup}\ Tr(\gamma A),
\end{equation}
where $\mathbb{S}^{+}_{d}$ denotes the collection of nonnegative elements in $\mathbb{S}_{d}$.
If there exists some $\beta >0$ such that $G(A)-G(B)\geq \beta tr[A-B]$ for any $A\geq B$, we call the G-normal
distribution non-degenerate, which is the case we consider throughout this paper.
\end{remark}

\begin{definition}
Let $\Omega=C_{0}^{d}([0,T])$, i.e., the space of all $\mathbb{R}^{d}$-valued continuous paths $(\omega_{t})_{t\in [0,T]}$ with $\omega_{0}=0$. The corresponding canonical process is $B_{t}(\omega)=\omega_{t}$, $t\in [0,T]$. $P_{0}$ is wiener measure. $\mathbb{F}=\{\mathcal{F}_{t}^{B}\}_{t\geq 0}$ is the filtration generated by $B$.  We let $\mathcal{H}:=L_{ip}(\Omega_{T})$ to be a linear space of random variables for each fixed $T\geq 0$, where $L_{ip}(\Omega_{T}):=\{\varphi(B_{t_{1}},\cdots,B_{t_{n}}):  n\geq 1, t_{1},\cdots,t_{n} \in [0,T],  \varphi \in C_{b,Lip}(\mathbb{R}^{d\times n})\}$.\par
$\left(\rm{i}\right)$ The G-expectation $\hat{E}$ is a sublinear expectation defined by
\begin{equation*}
 \hat{E}[X]:=\tilde{E}[\varphi(\sqrt{t_{1}-t_{0}}\xi_{1},\cdots,\sqrt{t_{n}-t_{n-1}}\xi_{n})],
\end{equation*}
for each $X= \varphi(B_{t_{1}}-B_{t_{0}}, B_{t_{2}}-B_{t_{1}},\cdots,B_{t_{n}}-B_{t_{n-1}})$,  where $(\xi_{i})_{i=1}^{n}$ are identically distributed d-dimensional G-normally distributed random vectors in a sublinear expectation space $(\tilde{\Omega},\tilde{\mathcal{H}},\tilde{E})$ such that $\xi_{i+1}$ is independent from $(\xi_{1},\cdots,\xi_{i})$  for each $i=1,2,\cdots,n-1$. $(\Omega,\mathcal{H},\hat{E})$ is called G-expectation space and the canonical process $\{B_{t}\}_{t \in [0,T]}$ in the sublinear space $(\Omega,\mathcal{H},\hat{E})$ is called a G-Brownian motion.

$\left(\rm{ii}\right)$ The conditional G-expectation $\hat{E}_{t}$ of $X \in L_{ip}(\Omega_{T})$ is defined by
\begin{equation*}
\hat{E}_{t_{j}}[\varphi(B_{t_{1}}-B_{t_{0}},B_{t_{2}}-B_{t_{1}},\cdots,B_{t_{n}}-B_{t_{n-1}})]:=\psi(B _{t_{1}}-B _{t_{0}},\cdots,B _{t_{j}}-B _{t_{j-1}}),
\end{equation*}
where $\psi(x_{1},\cdots,x_{j})=\tilde{E}[\varphi(x_{1},\cdots,x_{j},\sqrt{t_{j+1}-t_{j}}\xi_{j+1},\cdots,\sqrt{t_{n}-t_{n-1}}\xi_{n})]$.
\end{definition}
We denote by $L_{G}^{p}(\Omega_{T})$, $p\geq 1$, the completion of G-expectation space $L_{ip}(\Omega_{T})$ under the norm $\parallel X \parallel _{p,G}:=(\hat{E}[|X|^{p}])^{\frac{1}{p}}$.  For all $t\in [0,T]$, $\hat{E}[\cdot]$ and $\hat{E}_{t}[\cdot]$ are continuous mapping on $L_{ip}(\Omega_{T})$ endowed with the norm $\parallel \cdot \parallel_{1,G}$. Therefore, it can be extended continuously to $L_{G}^{p}(\Omega_{T})$.

\begin{definition}
A process $\{M_{t}\}_{t\geq 0}$ is called a G-martingale  if for each $t\in [0,T]$,
$M_{t}\in L_{G}^{1}(\Omega_{t})$ and for each $s\in [0,t]$, we have $\hat{E}_{s}[M_{t}]=M_{s}$.
\end{definition}
Now we introduce the It\^{o} integral and quadratic variation process with respect to G-Brownian motion in G-expectation space.

\begin{definition}
Let $p\geq1$ be fixed. For a given partition $\pi_{T}=\{t_{0},\cdots,t_{N}\}$ of $[0,T]$, we denote $M^{p,0}_{G}(0,T)$ as the collection of following type of simple processes:
\begin{equation*}
\eta_{t}(\omega)=\underset{k=0}{\overset{N-1}{\sum}}\xi_{k}(\omega)1_{[t_{k},t_{k+1})}(t),
\end{equation*}
where $\xi_{k}\in L_{ip}(\Omega_{t_{k}})$, $k=0,1,2,\cdots,N-1$. We denote by $M_{G}^{p}(0,T)$ the completion of $M^{p,0}_{G}(0,T)$ under the norm $\| \cdot \|_{M_{G}^{p}(0,T)}:= \{\hat{E}[\int^{T}_{0}|\cdot|^{p}dt]\}^{\frac{1}{p}}$.
\end{definition}

\begin{definition}
For each $\eta \in M_{G}^{2,0}(0,T)$, we define
\begin{equation*}
I(\eta)=\int^{T}_{0}\eta_{t}d B_{t}:=\underset{j=0}{\overset{N-1}{\sum}}\xi_{j}(B_{t_{j+1}}-B_{t_{j}}).
\end{equation*}
The mapping $I: M_{G}^{2,0}(0,T)\rightarrow L_{G}^{2}(\Omega_{T})$ is continuous and thus can be continuously extended to $M_{G}^{2}(0,T)$.
\end{definition}

\begin{definition}
The quadratic variation process of G-Brownian motion is defined by
\begin{equation*}
\langle B \rangle_{t}:= B_{t}^{2}-2\int^{t}_{0}B_{s}dB_{s},
\end{equation*}
which is a continuous, nondecreasing process.
\end{definition}

\begin{definition}
We now define the integral of a process $\eta\in M_{G}^{1}(0,T)$ with respect to $\langle B\rangle$ as following:
\begin{equation*}
Q_{0,T}(\eta)=\int^{T}_{0}\eta_{t}d\langle B \rangle_{t}:=\underset{j=0}{\overset{N-1}{\sum}}\xi_{j}(\langle B \rangle_{t_{j+1}}-\langle B \rangle_{t_{j}}):M_{G}^{1,0}(0,T)\rightarrow L_{G}^{1}(\Omega_{T}).
\end{equation*}
The mapping is continuous and can be extended to $M_{G}^{1}(0,T)$ uniquely.
\end{definition}
\par
Then we detail some results about the quasi-analysis theory constructed in \cite{r5}.
\begin{theorem}
There exists a weakly compact family $\mathcal{P}\subset \mathcal{M}_{1}(\Omega_{T})$, the collection of probability measures defined on $(\Omega_{T},\mathcal{B}(\Omega_{T}))$,  such that
\begin{equation*}
\hat{E}[X]=\underset{P\in\mathcal{P}}{sup}\ E_{p}[X],\  \forall X \in L_{ip}(\Omega_{T}),
\end{equation*}
$\mathcal{P}$ is called a  set of probability measures that represents $\hat{E}$.
\end{theorem}

\begin{definition}
We define the  capacity associated to $\mathcal{P}$, which is a weakly compact family of probability measure represents $\hat{E}$, as follow:
\begin{equation*}
\hat{c}(A):=\underset{P\in \mathcal{P}}{sup}\ P(A),\ A\in \mathcal{B}(\Omega_{T}),
\end{equation*}
$\hat{c}$ is also called the capacity induced by $\hat{E}$.
\end{definition}

Let $(\Omega^{0},\mathcal{F}^{0}=\{\mathcal{F}^{0}_{t}\},\mathcal{F},P^{0})$ be a filtered probability space, and $\{W_{t}\}$ be a d-dimensional Brownian motion under $P^{0}$. \cite{r5} proved that
$\mathcal{P}_{M}:=\{ P_{0} \circ X^{-1}|X_{t}=\int^{t}_{0}h_{s}dW_{s}, h\in L^{2}_{\mathcal{F}^{0}}([0,T];\Gamma^{\frac{1}{2}})\}$ represents G-expectation $\hat{E}$, where $\Gamma^{\frac{1}{2}}:=\{\gamma^{\frac{1}{2}}|\gamma \in \Gamma \}$ and $\Gamma$ is the set in the representation of $G(\cdot)$ of the formula $\eqref{11}$.

\begin{definition}
$\left(\rm{i}\right)$ Let $\hat{c}$ be the capacity induced by $\hat{E}$. A set $A\subset\Omega$ is polar if $\hat{c}(A)=0$. A property holds "quasi-surely"(q.s. for short) if it holds outside a polar set.

$\left(\rm{ii}\right)$ Let $X$ and $Y$ be two random variables, we say that $X$ is a version of $Y$, if $X=Y$ q.s.
\end{definition}
let $\parallel\psi \parallel_{p,G}=[\hat{E}(|\psi|^{p})]^{\frac{1}{p}}$ for $\psi \in C_{b}(\Omega_{T})$. The completion of $C_{b}(\Omega_{T})$ and $L_{ip}(\Omega_{T})$ under $\parallel \cdot \parallel_{p,G}$ are the same and we denote them by $L_{G}^{p}(\Omega_{T})$.

\subsection{Forward and Backward Stochastic Differential Equations Driven by G-Brownian Motion}

We consider the following stochastic differential equations driven by $d$ dimensional G-Brownian motion (G-SDE):
\begin{equation}\label{formula4}
X_{t}=X_{0}+\int^{t}_{0}b(s,X_{s})ds+\underset{i,j=1}{\overset{d}{\sum}}\int^{t}_{0}h_{ij}(s,X_{s})d\langle B^{i},B^{j}\rangle_{s}+\underset{j=1}{\overset{d}{\sum}}\int^{t}_{0}\sigma_{j}(s,X_{s})dB_{s}^{j},
\end{equation}
where $t\in[0,T]$, the initial condition $X_{0}\in\mathbb{R}^{n}$ is a given constant,  $b,h_{ij},\sigma_{j}$ are given functions satisfying $b(\cdot,x)$, $h_{ij}(\cdot,x)$, $\sigma_{j}(\cdot,x)\in {M}_{G}^{2}(0,T;\mathbb{R}^{n})$ for each $x\in\mathbb{R}^{n}$ and the Lipschitz condition, i.e., $|\phi(t,x)-\phi(t,x')|\leq K|x-x'|$, for each $t\in[0,T]$, $x$, $x'\in\mathbb{R}^{n}$, $\phi=b$, $h_{ij}$ and $\sigma_{j}$, respectively.  The solution is a process $X\in {M}^{2}_{G}(0,T;\mathbb{R}^{n})$ satisfying the G-SDE \eqref{formula4}.
\begin{theorem}(\cite{r3})
There exists a unique solution $X\in {M}^{2}_{G}(0,T;\mathbb{R}^{n})$ of the stochastic differential equation \eqref{formula4}.
\end{theorem}
Now we give the results about BSDEs driven by G-Brownian motion in the G-expectation space $(\Omega_{T},L_{G}^{1}(\Omega_{T}),\hat{E})$ with $\Omega_{T}=C_{0}([0,T],\mathbb{R}^{d})$ and $\bar{\sigma}^{2}=\hat{E}[B_{1}^{2}]\geq-\hat{E}[-B_{1}^{2}]=\underline{\sigma}^{2}>0$.
We consider the following type of G-BSDEs (we always use Einstein convention),
\begin{equation}\label{formula1}
Y_{t}=\xi+\int^{T}_{t}f(s,Y_{s},Z_{s})ds+\underset{i,j=1}{\overset{d}{\sum}}\int^{T}_{t}g_{ij}(s,Y_{s},Z_{s})d\langle B^{i},B^{j} \rangle_{s}-\int^{T}_{t}Z_{s}dB_{s}-(K_{T}-K_{t}),
\end{equation}
where $f(t,\omega,y,z)$, $g_{ij}(t,\omega,y,z):[0,T]\times\Omega_{T}\times\mathbb{R}\times \mathbb{R}^{d}\rightarrow\mathbb{R}$ satisfy the following properties: There exists some $\beta>1$ such that

$\rm{(H1)}$ for any y, z, $f(\cdot,\cdot,y,z$), $g_{ij}(\cdot,\cdot,y,z)\in M_{G}^{\beta}(0,T)$;

$\rm{(H2)}$ for some $L>0$,
\begin{equation*}
|f(t,\omega,y,z)-f(t,\omega,y',z')|+\underset{i,j=1}{\overset{d}{\sum}}|g_{ij}(t,\omega,y,z)-g_{ij}(t,\omega,y',z')|\leq L(|y-y'|+|z-z'|).
\end{equation*}

For simplicity, we denote by $\mathfrak{G}(0,T)$ the collection of processes $(Y,Z,K)$ such that $Y\in S_{G}^{\alpha}(0,T)$, $Z\in H_{G}^{\alpha}(0,T)$, $K$ is a decreasing G-martingale with $K_{0}=0$ and $K_{T}\in L_{G}^{\alpha}(\Omega_{T})$. Here $S_{G}^{\alpha}(0,T)$ is the completion of $S_{G}^{0}(0,T)=\{h(t,B_{t_{1}\wedge t,\cdots,B_{t_{n}\wedge t}}):t_{1},\cdots,t_{n}\in[0,T],h\in C_{b,lip(\mathbb{R}^{n+1})}\}$ under $\parallel\cdot\parallel_{s_{G}^{p}}=\{\hat{E}[sup_{t\in[0,T]}|\eta_{t}|^{p}]\}^{\frac{1}{p}}$ and $H_{G}^{p}(0,T)$ is the completion of $M_{G}^{0}(0,T)$ under $\|\cdot\|_{H_{G}^{P}}=\{\hat{E}[(\int^{T}_{0}|\eta_{s}|^{2}ds)^{p/2}]\}^{1/p}$.
\begin{definition}
Let $\xi\in L_{G}^{\beta}(\Omega_{T})$ with $\beta> 1$, $f$ and $g_{ij}$ satisfy $(H1)$ and $(H2)$. A triplet of processes $(Y,Z,K)$ is called a solution of equation \eqref{formula1} if for some $1<\alpha\leq \beta$ the following properties hold:

$\rm{(a)}$ $(Y,Z,K)\in\mathfrak{G}_{G}^{\alpha}(0,T)$;

$\rm{(b)}$
\begin{equation*}
Y_{t}=\xi+\int^{T}_{t}f(s,Y_{s},Z_{s})ds+\underset{i,j=1}{\overset{d}{\sum}}\int^{T}_{t}g_{ij}(s,Y_{s},Z_{s})d\langle B^{i},B^{j} \rangle_{s}-\int^{T}_{t}Z_{s}dB_{s}-(K_{T}-K_{t}).
\end{equation*}
\end{definition}\label{theorem1}
\begin{theorem}(\cite{r20})
Assume that $\xi \in L_{G}^{\beta}(\Omega_{T})$ and $f, g$ satisfy $(H1)$ and $(H2)$ for some $\beta > 1$. Then equation \eqref{formula1} has a unique solution $(Y,Z,K)$. Moreover,  for any $1 < \alpha < \beta$ we have $Y \in S_{G}^{\alpha}(0,T)$, $Z \in H_{G}^{\alpha}(0,T; \mathbb{R}^{d})$ and $K_{T}\in L_{G}^{\alpha}(\Omega_{T})$.
\end{theorem}

We have the following estimates.
\begin{proposition}(\cite{r20})
Let $\xi \in L_{G}^{\beta}(\Omega_{T})$ and $f$, $g_{ij}$ satisfy $(H1)$ and $(H2)$ for some $\beta > 1$. Assume that $(Y,Z,K) \in \mathfrak{G}_{G}^{\alpha}(0,T)$ for some $1 < \alpha < \beta$ is a solution of equation \eqref{formula1}. Then
\par
$\rm{(i)}$ There exists a constant $C_{\alpha}:= C(\alpha,T,G,L)> 0$ such that
\begin{equation*}
|Y_{t}|^{\alpha}\leq C_{\alpha} \hat{E}_{t}[|\xi|^{\alpha}+\int^{T}_{t}|h_{s}^{0}|^{\alpha}ds],
\end{equation*}
\begin{equation*}
\hat{E}[(\int^{T}_{0}|Z_{s}|^{2}ds)^{\frac{\alpha}{2}}]\leq C_{\alpha}\{\hat{E}[\underset{t \in[0,T]}{\sup}|Y_{t}|^{\alpha}]+(\hat{E}[\underset{t \in [0,T]}{\sup}|Y_{t}|^{\alpha}])^{\frac{1}{2}}(\hat{E}[(\int_{0}^{T}h_{s}^{0}ds)^{\alpha}])^{\frac{1}{2}}\},
\end{equation*}
\begin{equation*}
\hat{E}[|K_{T}|^{\alpha}]\leq C_{\alpha}\{\hat{E}[\underset{t \in[0,T]}{\sup}|Y_{t}|^{\alpha}]+\hat{E}[(\int_{0}^{T}h_{s}^{0}ds)^{\alpha}]\},
\end{equation*}
where $h_{s}^{0}=|f(s,0,0)|+\underset{i,j=1}{\overset{d}{\sum}}|g_{ij}(s,0,0)|$.
\par
$\rm{(ii)}$ For any given $\alpha < \alpha' < \beta$, there exists a constant $C_{\alpha,\alpha'}$ depending on $\alpha,\alpha',T,G,L$ such that
\begin{align*}
\hat{E}[\underset{t\in[0,T]}{\sup}|Y_{t}|^{\alpha}]&\leq C_{\alpha,\alpha'}\{\hat{E}[\underset{t\in[0,T]}{\sup}\hat{E}_{t}[|\xi|^{\alpha}]]\\
 &~~~~+(\hat{E}[\underset{t\in[0,T]}{\sup}\hat{E}_{t}[(\int_{0}^{T}h_{s}^{0}ds)^{\alpha'}]])^{\frac{\alpha}{\alpha'}}
 +\hat{E}[\underset{t\in[0,T]}{\sup}\hat{E}_{t}[(\int_{0}^{T}h_{s}^{0}ds)^{\alpha'}]]\}.
\end{align*}
\end{proposition}
\begin{proposition}(\cite{r21})
Let $\xi^{i} \in L_{G}^{\beta}(\Omega_{T})$, $i=1,2$, and $f^{i}$, $g_{ij}^{i}$ satisfy $(H1)$ and $(H2)$ for some $\beta > 1$. Assume that $(Y^{i},Z^{i},K^{i})\in \mathfrak{G}_{G}^{\alpha}(0,T)$ for some $1 < \alpha < \beta$ are the solutions of \eqref{formula1} corresponding to $\xi^{i}$, $f^{i}$ and $g_{ij}^{i}$. Set $\hat{Y}_{t}=Y_{t}^{1}-Y_{t}^{2}$, $\hat{Z}_{t}=Z_{t}^{1}-Z_{t}^{2}$ and $\hat{K}_{t}=K_{t}^{1}-K_{t}^{2}$. Then
\par
$\rm{(i)}$ There exists a constant $C_{\alpha}:=C(\alpha,T,G,L)> 0$ such that
\begin{equation*}
|\hat{Y}_{t}|^{\alpha}\leq C_{\alpha}\hat{E}_{t}[|\hat{\xi}|^{\alpha}+\int_{t}^{T}|\hat{h}_{s}|^{\alpha}ds],
\end{equation*}
where $\hat{\xi}=\xi^{1}-\xi^{2}$, $\hat{h}_{s}=|f^{1}(s,Y_{s}^{2},Z_{s}^{2})-f^{2}(s,Y_{s}^{2},Z_{s}^{2})|
+\underset{i,j=1}{\overset{d}{\sum}}|g_{ij}^{1}(s,Y_{s}^{2},Z_{s}^{2})-g_{ij}^{2}(s,Y_{s}^{2},Z_{s}^{2})|$.
\par
$\rm{(ii)}$ For any given $\alpha'$ with $\alpha < \alpha' < \beta$, there exists a constant $C_{\alpha,\alpha'}$ depending on $\alpha,\alpha',T,G,L$ such that
\begin{align*}
\hat{E}[\underset{t \in [0,T]}{\sup}|\hat{Y}_{t}|^{\alpha}]&\leq C_{\alpha,\alpha'}\{\hat{E}[\underset{t\in [0,T]}{\sup}\hat{E}_{t}[|\hat{\xi}|^{\alpha}]]\\
&~~~~+(\hat{E}[\underset{t\in [0,T]}{\sup}\hat{E}_{t}[(\int_{0}^{T}\hat{h}_{s}ds)^{\alpha'}]])^{\frac{\alpha}{\alpha'}}+\hat{E}[\underset{t\in [0,T]}{\sup}\hat{E}_{t}[(\int_{0}^{T}\hat{h}_{s}ds)^{\alpha'}]]\}.
\end{align*}
\end{proposition}
\par

\begin{theorem}(\cite{r21})
Let $(Y_{t}^{i},Z_{t}^{i},K_{t}^{i})_{t\leq T}$, $i=1,2$, be the solutions of the following G-BSDEs:
\begin{equation*}
Y_{t}^{i}=\xi^{i}+\int_{t}^{T}f_{i}(s,Y_{s}^{i},Z_{s}^{i})ds+\int_{t}^{T}g_{i}(s,Y_{s}^{i},Z_{s}^{i})d \langle B \rangle_{s}-\int_{t}^{T}Z_{s}^{i}dB_{s}-(K_{T}^{i}-K_{t}^{i}),
\end{equation*}
where $\xi^{i} \in L_{G}^{\beta}(\Omega_{T})$, $f_{i}$, $g_{i}$ satisfy $(H1)$ and $(H2)$ with $\beta > 1$. If $\xi^{1}\geq \xi^{2}$, $f_{1}\geq f_{2}$, $g_{1}\geq g_{2}$, then $Y_{t}^{1}\geq Y_{t}^{2}$.
\end{theorem}
\begin{theorem}(\cite{r21})
Let $(Y_{t}^{i},Z_{t}^{i},K_{t}^{i})_{t\leq T}$, $i=1,2$, be the solutions of the following G-BSDEs:
\begin{equation*}
Y_{t}^{i}=\xi^{i}+\int_{t}^{T}f_{i}(s,Y_{s}^{i},Z_{s}^{i})ds+\int_{t}^{T}g_{i}(s,Y_{s}^{i},Z_{s}^{i})d \langle B \rangle_{s}-\int_{t}^{T}Z_{s}^{i}dB_{s}-(K_{T}^{i}-K_{t}^{i})+V^{i}_{T}-V_{t}^{i},
\end{equation*}where $\xi^{i} \in L_{G}^{\beta}(\Omega_{T})$, $f_{i}$, $g_{i}$ satisfy $(H1)$ and $(H2)$, $(V_{t}^{i})_{t\leq T}$ are RCLL processes such that $\hat{E}[\underset{t\in[0,T]}{\sup}|V_{t}^{i}|^{\beta}] < \infty$  with $\beta > 1$. If $\xi^{1}\geq \xi^{2}$, $f_{1}\geq f_{2}$, $g_{1}\geq g_{2}$, $V_{t}^{1}-V_{t}^{2}$ is an increasing process, then $Y_{t}^{1}\geq Y_{t}^{2}$.
\end{theorem}

\section{A DPP for Stochastic Optimal Control Problems under G-Expectation}
Now we introduce the setting for stochastic optimal control problems under $G$-expectation. We suppose that the control state space $V$ is a compact metric space. Let the set of admissible control processes $\mathcal{U}$ for the player be a set of $V$-valued stochastic processes in $M_{G}^{\beta}([t,T];\mathbb{R}^{n})$ with $\beta > 2$ and $t\in[0,T]$. For a given admissible control $\upsilon(\cdot)\in\mathcal{U}$, the corresponding orbit which regards $t$ as the initial time and $\xi\in L_{G}^{2}(\Omega_{t};\mathbb{R}^{n})$ as the initial state, is defined by the solution of the following type of G-SDE:

\begin{equation}\label{formula3}
\left \{
\begin{array}{lr}
dX_{s}^{t,\xi;\upsilon}=b(s,X_{s}^{t,\xi;\upsilon},\upsilon_{s})ds+\underset{i,j=1}
{\overset{d}{\sum}}h_{ij}(s,X_{s}^{t,\xi;\upsilon},\upsilon_{s})d\langle B^{i},B^{j}\rangle_{s}+\underset{j=1}{\overset{d}{\sum}}\sigma_{j}(s,X_{s}^{t,\xi;\upsilon},\upsilon_{s})dB_{s}^{j},\\
s\in[t,T], \\
X_{t}^{t,\xi;\upsilon}=\xi, \\
\end{array}
\right.
\end{equation}
where $b$, $h_{ij}$, $\sigma_{j}:[0,T]\times \mathbb{R}^{n}\times \mathcal{U}\rightarrow \mathbb{R}^{n}$ are deterministic functions and satisfy the following conditions $\rm{(H3)}$:
\par
$\rm{(A1)}$ $h_{ij}=h_{ji}$ for $1 \leq i,j \leq d$;\par
$\rm{(A2)}$ For every fixed $(x,\upsilon)\in\mathbb{R}^{n}\times \mathcal{U}$, $b(\cdot,x,\upsilon)$, $h_{ij}(\cdot,x,\upsilon)$, $\sigma_{j}(\cdot,x,\upsilon)$ are continuous in $t$;\par
$\rm{(A3)}$ There exists a constant $L > 0$, for any $t\in[0,T]$, $x$, $x'\in \mathbb{R}^{n}$, $\upsilon$, $\upsilon' \in \mathcal{U}$ such that
\begin{align*}
&|b(t,x,\upsilon)-b(t,x',\upsilon')|+\underset{i,j=1}{\overset{d}{\Sigma}}
|h_{ij}(t,x,\upsilon)-h_{ij}(t,x',\upsilon')|+\underset{j=1}{\overset{d}{\Sigma}}|\sigma_{j}(t,x,\upsilon)-\sigma_{j}(t,x',\upsilon')|\\
&\leq L(|x-x'|+|\upsilon-\upsilon'|).
\end{align*}
\par
From the assumption $\rm{(H3)}$, we can get global linear growth conditions for $b$, $h_{ij}$, $\sigma_{j}$, i.e., there exists  $C> 0$ such that, for $t\in[0,T]$, $x\in\mathbb{R}^{n}$, $\upsilon\in \mathcal{U}$, $|b(t,x,\upsilon)|+\underset{i,j=1}{\overset{d}{\sum}}|h_{ij}(t,x,\upsilon)|+\underset{j=1}{\overset{d}{\sum}}|\sigma_{j}(t,x,\upsilon)|\leq C(1+|x|+|\upsilon|)$. Obviously, under the above assumptions, for any $\upsilon(\cdot)\in\mathcal{U}$, G-SDE \eqref{formula3} has a unique solution. Moreover, we have the following estimates :
\begin{proposition}
Let $\xi$, $\xi'\in L_{G}^{p}(\Omega_{t};\mathbb{R}^{n})$ with $p\geq 2$, $\upsilon(\cdot)$, $\upsilon'(\cdot) \in \mathcal{U}$, $t \in [0,T]$ , $\delta \in [0,T-t]$, then we have
\begin{equation*}
\hat{E}_{t}[|X_{t+\delta}^{t,\xi;\upsilon}-X_{t+\delta}^{t,\xi';\upsilon'}|^{p}]\leq C(|\xi-\xi'|^{p}+\int_{t}^{t+\delta}\hat{E}_{t}|\upsilon_{r}-\upsilon_{r}'|^{p}dr),
\end{equation*}
\begin{equation*}
\hat{E}_{t}[|X_{t+\delta}^{t,\xi;\upsilon}|^{p}]\leq C(1+|\xi|^{p}),
\end{equation*}
\begin{equation*}
\hat{E}_{t}[\underset{s\in [t,t+\delta]}{\sup}|X_{s}^{t,\xi;\upsilon}-\xi|^{p}]\leq C(1+|\xi|^{p})\delta^{\frac{p}{2}},
\end{equation*}
where $C$ depends on $L,G,p,n,T$.
\end{proposition}
\proof The proof is similar to the proof of Proposition 4.1 in \cite{r21}.
\endproof
\par
Now we give bounded functions $\Phi:\mathbb{R}^{n}\rightarrow\mathbb{R}$, $f:[0,T]\times \mathbb{R}^{n}\times\mathbb{R}\times\mathbb{R}^{d}\times\mathcal{U}\rightarrow\mathbb{R}$, $g_{ij}:[0,T]\times \mathbb{R}^{n}\times\mathbb{R}\times\mathbb{R}^{d}\times\mathcal{U}\rightarrow\mathbb{R}$ satisfy the following conditions:
$\rm{(H4)}$ \par
$\rm{(i)}$  $g_{ij}=g_{ji}$ for $1 \leq i,j \leq d$.\par
$\rm{(ii)}$ For every fixed $(x,y,z,\upsilon)\in \mathbb{R}^{n}\times \mathbb{R}\times \mathbb{R}^{n}\times \mathcal{U}$, $f(\cdot,x,y,z,\upsilon)$ and $g_{ij}(\cdot,x,y,z,\upsilon)$ are continuous in $t$, $1 \leq i,j \leq d$.\par
$\rm{(iii)}$ There exist a constant $ L > 0$, for $t\in[0,T]$, $x$, $x'\in \mathbb{R}^{n}$, $y$, $y'\in \mathbb{R}$, $z$, $z'\in \mathbb{R}^{d}$, $\upsilon$, $\upsilon'\in \mathcal{U}$, such that
\begin{equation*}
|\Phi(x)-\Phi(x')| \leq L (|x-x'|),
\end{equation*}
\begin{align*}
&|f(t,x,y,z,\upsilon)-f(t,x',y',z',\upsilon')|+\underset{i,j=1}{\overset{d}{\Sigma}}|g_{ij}(t,x,y,z,\upsilon)-g_{ij}(t,x',y',z',\upsilon')|\\
&\leq L (|x-x'|+|y-y'|+|z-z'|+|\upsilon-\upsilon'|).
\end{align*}

From $\rm{(H4)}$, we have that $\Phi$, $f$ and $g_{ij}$ also satisfy global linear growth condition in $x$, i.e., there exists $C > 0$, such that for all $0 \leq t \leq T$, $\upsilon \in \mathcal{U}$, $x \in \mathbb{R}^{n}$,
\begin{equation*}
|\Phi(x)|+|f(t,x,0,0,\upsilon)|+|g_{ij}(t,x,0,0,\upsilon)|\leq C(1+|x|+|\upsilon|).
\end{equation*}
For any $\upsilon \in \mathcal{U}$ and $\xi \in L_{G}^{2}(\Omega_{t},\mathbb{R}^{n})$, the mappings $f(s,x,y,z,\upsilon):=f(s,X_{s}^{t,\xi;\upsilon},y,z,\upsilon_{s})$ and $g_{ij}(s,x,y,z,\upsilon)=g_{ij}(s,X_{s}^{t,\xi;\upsilon},y,z,\upsilon_{s})$, where $(s,y,z) \in [0,T]\times \mathbb{R}\times \mathbb{R}^{d}$, satisfy the conditions of Theorem 2.16 on the interval $[t,T]$. Therefore, there exists a unique solution for the following G-BSDE:
\begin{align}\label{formula6}
Y_{s}^{t,\xi;\upsilon}=&\Phi(X_{T}^{t,\xi;\upsilon})
+\int_{s}^{T}f(r,X_{r}^{t,\xi;\upsilon},Y_{r}^{t,\xi;\upsilon},Z_{r}^{t,\xi;\upsilon},\upsilon_{r})dr
-\int^{T}_{s}Z_{r}^{t,\xi;\upsilon}dB_{r}-(K_{T}^{t,\xi;\upsilon}-K_{s}^{t,\xi;\upsilon})\nonumber\\
&+\underset{i,j=1}{\overset{d}{\sum}}\int^{T}_{s}g_{ij}(r,X_{r}^{t,\xi;\upsilon},
Y_{r}^{t,\xi;\upsilon},Z_{r}^{t,\xi;\upsilon},\upsilon_{r})d\langle B^{i},B^{j} \rangle_{r},
\end{align}
where $X^{t,\xi;\upsilon}$ is introduced by \eqref{formula3}.
\begin{proposition}
 For each $ \xi$, $\xi' \in L_{G}^{p}(\Omega_{t};\mathbb{R}^{n})$ with $p\geq 2$ and $\upsilon(\cdot)$, $\upsilon'(\cdot)\in \mathcal{U}$ we have
\begin{equation*}
|Y_{t}^{t,\xi;\upsilon}-Y_{t}^{t,\xi';\upsilon}|\leq C |\xi-\xi'| ,
\end{equation*}
\begin{equation*}
|Y_{t}^{t,\xi;\upsilon}|\leq C(1+|\xi|),
\end{equation*}
\begin{equation*}
|Y_{t}^{t,\xi;\upsilon}-Y_{t}^{t,\xi;\upsilon'}|\leq C(\int_{t}^{T}\hat{E}_{t}|\upsilon(r)-\upsilon'(r)|^{2}dr)^{\frac{1}{2}},
\end{equation*}
where $C$ depends on $L$, $G$, $n$ and $T$.
\end{proposition}
\proof The proof is similar to the Proposition 4.2 in \cite{r21}.
\endproof\par
Given a control process $\upsilon(\cdot) \in \mathcal{U}$, we introduce an associated cost functional
\begin{equation*}
J(t,x;\upsilon)=Y_{t}^{t,x;\upsilon},\ (t,x)\in [0,T]\times \mathbb{R}^{n},
\end{equation*}
where the process $Y^{t,\xi;\upsilon}_{t}$ is defined by G-BSDE \eqref{formula6}.
Similar to the proof of Theorem 4.4 in \cite{r21}, we have that for $t \in [0,T]$,  $\xi \in L_{G}^{2}(\Omega_{t},\mathbb{R}^{n})$,
\begin{equation*}
J(t, \xi; \upsilon):=Y_{t}^{t, \xi;\upsilon}.
\end{equation*}
But we are more interest in the case when $\xi= x$.

Now we  define the value function as follow:
\begin{equation}\label{18}
u(t,x):=\underset{\upsilon(\cdot)\in \mathcal{U}}{\sup}\ J(t,x;\upsilon).
\end{equation}
\begin{proposition}
$u(t,x)$ is a deterministic function of $(t,x)$.
\end{proposition}
\proof For a partition of $[t,s]$: $t=t_{0}<t_{1}<\cdots<t_{N}=s$, $p \geq 2$, $t\leq s \leq T$, we denote  $L_{ip}(\Omega^{t}_{s}):=\{\varphi(B_{t_{1}}-B_{t},\cdots,B_{t_{n}}-B_{t}):n\geq 1,t_{1},\cdots,t_{n}\in[t,s],\varphi \in C_{b,Lip}(\mathbb{R}^{d\times n})\}$, $M_{G}^{p,0,t}(t,s;\mathbb{R}^{n})$ by the collection of simple processes $\eta(r)=\underset{k=0}{\overset{N-1}{\sum}}\xi_{k}1_{[t_{k},t_{k+1})}(r)$, where $\xi_{k}\in L_{ip}(\Omega^{t}_{t_{k}};\mathbb{R}^{n})$, $k=0,1,2,\cdots,N-1$ and $M_{G}^{p,t}(t,s;\mathbb{R}^{n})$ by the completion of $M^{p,0,t}_{G}(t,s;\mathbb{R}^{n})$ under the norm $\| \eta\|_{M_{G}^{p}(t,s;\mathbb{R}^{n})}:= \{\hat{E}[\int^{s}_{t}|\eta(r)|^{p}dr]\}^{\frac{1}{p}}$. Use the similar method in Lemma 43 of \cite{r5}, we can prove for $\upsilon \in M_{G}^{p}(t,s;\mathbb{R}^{n})$ is a $V$-valued process, there exists $\{u=\underset{i=1}{\overset{N}{\Sigma}}1_{A_{i}}u^{i}\}_{N\in \mathbb{N}}$,  $u^{i}\in M_{G}^{p,t}(t,s;\mathbb{R}^{n})$ is a $V$-valued process, $A_{i}$ is a partition of $\mathcal{B}(\Omega_{t})$ such that $u\rightarrow v$ under the norm $\| \eta\|_{M_{G}^{p}(t,s;\mathbb{R}^{n})}$, $N\rightarrow \infty$.  When $\upsilon(s)\in M_{G}^{p,t}(t,s;\mathbb{R}^{n})$, we note that $J(t,x;\upsilon)$ is a deterministic function of $(t,x)$, because $b$, $h_{ij}$, $\sigma_{j}$, $\Phi$, $f$ and $g_{ij}$ are deterministic functions and $\tilde{B}_{s}:=B_{t+s}-B_{t}$ is a G-Brownian motion. So we need to construct a sequence of admissible controls $\{\tilde{\upsilon}^{i}(\cdot)\}$ of the form
\begin{equation*}
\tilde{\upsilon}_{s}^{i}=\underset{j=1}{\overset{N_{i}}{\sum}} \upsilon_{s}^{ij}1_{A_{ij}}
\end{equation*}
satisfying $\underset{i \rightarrow \infty}{lim}J(t,x;\tilde{\upsilon}^{i}(\cdot))=u(t,x)$,  where $\upsilon^{ij}(\cdot)\in M_{G}^{p,t}(t,s;\mathbb{R}^{n})$ is a $V$-valued processes and $\{A_{ij}\}^{N_{i}}_{j=1}$ is a partition of $\mathcal{B}(\Omega_{t})$.
Firstly, there exists $\{\upsilon^{k}\}_{k \geq 1} \subset \mathcal{U}$, such that $u(t,x)= \underset{k \geq 1}{\sup}\ J(t,x;\upsilon^{k})$. Then we define $\upsilon,\upsilon' \in \mathcal{U}$,
\begin{equation*}
(\upsilon\vee \upsilon')_{s}=\left \{
\begin{array}{rl}
&0,s\in[0,t];\\
&\upsilon_{s},s\in(t,T],\ on\ \{J(t,x;\upsilon)\geq J(t,x;\upsilon')\};\\
&\upsilon_{s}',s\in(t,T],\ on\ \{J(t,x;\upsilon)< J(t,x;\upsilon')\}.\\
\end{array}
\right.
\end{equation*}
Therefore,
 $$J(t,x;\upsilon\vee \upsilon')\geq J(t,x;\upsilon) \vee J(t,x;\upsilon').$$
Set $\bar{\upsilon}^{1}:=\upsilon^{1}\vee \upsilon^{1}$, $\bar{\upsilon}^{k}:=\bar{\upsilon}^{k-1}\vee \upsilon^{i}$, $i \geq 2$. So $u(t,x)=\underset{k \rightarrow \infty}{lim}\ J(t,x;\bar{\upsilon}^{k})$. Without loss of generality, suppose
$\hat{E}[(u(t,x)-J(t,x;\bar{\upsilon}^{k}))^{2}] \leq 1/k,\ k \geq 1.$
 We denote
\begin{equation*}
\tilde{\upsilon}^{k}_{s}= \underset{j,k=0}{\overset{N_{i}-1}{\sum}} \bar{\upsilon}_{j,k}(s)1_{A_{j}^{k}},
\end{equation*}
where $\bar{\upsilon}_{j,k}\in M_{G}^{p,t}(t,s;\mathbb{R}^{n})$ is a $V$-values provess, $\{A_{j}^{k}\}_{0 \leq j \leq N_{k}-1}$ is a partition of $\mathcal{B}(\Omega_{t})$. Then we can suppose for $k \geq 1$, $\hat{E}[\int_{t}^{T} |\bar{\upsilon}_{s}^{k}-\tilde{\upsilon}_{s}^{k}|^{2}ds] \leq \frac{1}{Ck}$. From Proposition 3.2, we have
\begin{equation*}
\hat{E}[|J(t,x;\bar{\upsilon}^{k})-J(t,x;\tilde{\upsilon}^{k})|^{2}] \leq C \hat{E}[\int_{t}^{T}|\bar{\upsilon}^{k}_{s}-\tilde{\upsilon}^{k}_{s}|^{2}ds] \leq \frac{1}{k}.
\end{equation*}
Therefore, $\hat{E}[|u(t,x)-J(t,x;\tilde{\upsilon}^{k})|^{2}]\leq \frac{4}{k}$. Then we have
 $$J(t,x;\tilde{\upsilon}^{k})=\underset{j=0}{\overset{N_{k}-1}{\sum}}1_{A_{j}^{k}}J(t,x;\bar{\upsilon}_{j,k})\leq u(t,x).$$
  Now we suppose that
$$J(t,x,\tilde{\upsilon})\leq max_{0 \leq j \leq N_{k}-1}J(t,x;\bar{\upsilon}_{j,k})=J(t,x;\bar{\upsilon}_{j',k}).$$
 Because $\hat{E}[|J(t,x;\tilde{\upsilon}^{k})-u(t,x)|^{2}] \rightarrow 0$, we have
 $$u(t,x)=\underset{k \rightarrow \infty}{lim} J(t,x;\bar{\upsilon}_{j',k}),\ q.s..$$
  Hence $\hat{E}[u(t,x)]=u(t,x)$. We have finished the proof.

\endproof
\begin{lemma}
For any $t \in [0,T]$, $x$, $x'\in \mathbb{R}^{n}$, we have
\begin{equation}\label{12}
|u(t,x)-u(t,x')|\leq C |x-x'|,
\end{equation}
\begin{equation}\label{13}
|u(t,x)|\leq C(1+|x|).
\end{equation}
\end{lemma}
\proof
By Proposition 3.2, we have for $\upsilon(\cdot)\in \mathcal{U}$,
\begin{equation*}
|J(t,x;\upsilon(\cdot))|\leq C(1+|x|),
\end{equation*}
\begin{equation*}
|J(t,x;\upsilon(\cdot))-J(t,x';\upsilon(\cdot))|\leq C|x-x'|.
\end{equation*}
Then, $\forall \varepsilon >0$, there exist $\upsilon(\cdot)$, $\upsilon'(\cdot)\in \mathcal{U}$ such that
\begin{equation*}
J(t,x;\upsilon(\cdot))\leq u(t,x) \leq J(t,x;\upsilon(\cdot))+\varepsilon,
\end{equation*}
\begin{equation*}
J(t,x';\upsilon'(\cdot))\leq u(t,x') \leq J(t,x';\upsilon'(\cdot))+\varepsilon.
\end{equation*}
Now we have
\begin{equation*}
-C(1+|x|)\leq J(t,x;\upsilon(\cdot)) \leq u(t,x) \leq J(t,x;\upsilon(\cdot))+\varepsilon \leq C(1+|x|)+\varepsilon.
\end{equation*}
So we get \eqref{13}. Similarly, we obtain
\begin{equation*}
J(t,x;\upsilon'(\cdot))-J(t,x';\upsilon'(\cdot))-\varepsilon\leq u(t,x)-u(t,x') \leq J(t,x;\upsilon(\cdot)) -J(t,x';\upsilon(\cdot))+\varepsilon.
\end{equation*}
Then
\begin{equation*}
-C|x-x'|-\varepsilon \leq |u(t,x)-u(t,x')| \leq C|x-x'|+\varepsilon.
\end{equation*}
Thus we have proved \eqref{12}.
\endproof
\begin{lemma}
For any $t\in [0,T]$, $\zeta \in L_{G}^{2}(\Omega_{t};\mathbb{R}^{n})$ and $\zeta$ is $\mathcal{F}^{B}_{t}$ measurable, we have $\forall \upsilon(\cdot) \in \mathcal{U}$,
\begin{equation}\label{14}
 u(t,\zeta)\geq Y_{t}^{t,\zeta;\upsilon}.
\end{equation}
Conversely, $\forall \varepsilon  >0$, there exists a $\upsilon(\cdot) \in \mathcal{U}$, such that
\begin{equation}\label{3.7}
u(t,\zeta) \leq Y_{t}^{t,\zeta;\upsilon}+\varepsilon.
\end{equation}
\end{lemma}
\proof
We already know that $u(t,x)$ is continuous with respect to $x$ and $Y_{t}^{t,\zeta;\upsilon}$ is continuous with respect to $(\zeta,\upsilon(\cdot))$. We want to prove \eqref{14}, only need to discuss the simple random variables $\zeta$ of the form
\begin{equation*}
\zeta=\underset{i=1}{\overset{N}{\Sigma}} 1_{A_{i}}x_{i},
\end{equation*}
and $\upsilon(\cdot)$ of the form
\begin{equation*}
\upsilon(\cdot)=\underset{i=1}{\overset{N}{\Sigma}} 1_{A_{i}}\upsilon^{i}(\cdot).
\end{equation*}
Here $i=1,2,...,N$, $x_{i}\in \mathbb{R}^{n}$, $\upsilon^{i}\in M_{G}^{p,t}(t,s;\mathbb{R}^{n})$ and $\{A_{i}\}_{i=1}^{N}$ is a $\mathcal{B}(\Omega_{t})$-partition. Then from the same technique used in the proof of Theorem 4.4 in \cite{r21}, we have
\begin{equation*}
Y_{t}^{t,\zeta;\upsilon}=\underset{i=1}{\overset{N}{\Sigma}}1_{A_{i}}Y_{t}^{t,x^{i};\upsilon^{i}} \leq \underset{i=1}{\overset{N}{\Sigma}}1_{A_{i}}u(t,x_{i})=u(t,\underset{i=1}{\overset{N}{\Sigma}}1_{A_{i}}x_{i})=u(t,\zeta).
\end{equation*}
So we have proved \eqref{14}.
Now we prove \eqref{3.7} in a similar way. We first construct a random variable $\eta \in L_{G}^{2}(\Omega_{t};\mathbb{R}^{n})$,
\begin{equation*}
\eta=\underset{i=1}{\overset{N}{\Sigma}} x_{s}^{i} 1_{A_{i}},
\end{equation*}
where $(A_{i})_{i=1}^{N}$ is a $\mathcal{B}(\Omega_{t})$-partition and $x_{i}\in \mathbb{R}^{n}$,  such that $|\eta - \zeta| \leq \frac{\varepsilon}{3C}$. Then we have
\begin{equation*}
|Y_{t}^{t,\eta;\upsilon}-Y_{t}^{t,\zeta;\upsilon}| \leq \frac{\varepsilon}{3},
\end{equation*}
\begin{equation*}
|u(t,\zeta)-u(t,\eta)|\leq \frac{\varepsilon}{3},
\end{equation*}
for $\upsilon(\cdot) \in \mathcal{U}$.
Now, we chose a  control $\upsilon^{i}(\cdot)\in M_{G}^{p,t}(t,s;\mathbb{R}^{n})$, such that $u(t,x_{i}) \leq Y_{t}^{t,x^{i};\upsilon^{i}}+\frac{\varepsilon}{3}$. Set $\upsilon(\cdot):= \underset{i=1}{\overset{N}{\sum}}\upsilon^{i}(\cdot)1_{A_{i}}$. Finally, we get

\begin{align*}
Y_{t}^{t\zeta;\upsilon}&\geq -|Y_{t}^{t,\eta;\upsilon}-Y_{t}^{t,\zeta;\upsilon}|+Y_{t}^{t,\eta;\upsilon}\\
&\geq -\frac{\varepsilon}{3}+\underset{i=1}{\overset{N}{\sum}}Y_{t}^{t,x^{i};\upsilon^{i}}1_{A_{i}}\\
&\geq -\frac{\varepsilon}{3}+\underset{i=1}{\overset{N}{\sum}} (u(t,x_{i})-\frac{\varepsilon}{3})1_{A_{i}}\\
&=-\frac{2\varepsilon}{3}+\underset{i=1}{\overset{N}{\sum}}u(t,x_{i})1_{A_{i}}\\
&=-\frac{2\varepsilon}{3}+u(t,\eta)\geq -\varepsilon+ u(t,\zeta).
\end{align*}
So we have \eqref{3.7}.
\endproof
Now we give a type of DPP for our stochastic optimal control problems. Firstly, we define a family of backward semigroups associated with the G-BSDE \eqref{formula6}. Given the initial data $(t,x)$, a positive number $\delta \leq T-t$ and a random variable $\eta \in L_{G}^{p}(\Omega;\mathbb{R})$ with $p > 1$, we set
\begin{equation*}
G_{t,t+\delta}^{t,x;\upsilon}[\eta]:=Y_{s}^{t,x;\upsilon},
\end{equation*}
where $(Y_{s}^{t,x;\upsilon})_{t \leq s \leq t+\delta}$ is the solution of the following G-BSDE with the time horizon $t+\delta$:
\begin{align*}
Y_{s}^{t,x;\upsilon}&=\eta+\int_{s}^{t+\delta}f(r,X_{r}^{t,x;\upsilon},Y_{r}^{t,x;\upsilon},
Z_{r}^{t,x;\upsilon},\upsilon_{r})dr-\int^{t+\delta}_{s}Z_{r}^{t,x;\upsilon}dB_{r}-(K_{T}^{t,x;\upsilon}-K_{t}^{t,x;\upsilon})\\
&~~~~+\underset{i,j=1}{\overset{d}{\sum}}\int^{t+\delta}_{s}g_{ij}(r,X_{r}^{t,x;\upsilon},
Y_{r}^{t,x;\upsilon},Z_{r}^{t,x;\upsilon},\upsilon_{r})d\langle B^{i},B^{j} \rangle_{r}.
\end{align*}
Obviously, for the solution $Y^{t,x;\upsilon}_{\cdot}$ of G-BSDE \eqref{formula6}, we have
\begin{equation*}
G_{t,T}^{t,x;\upsilon}[\Phi(X_{T}^{t,x;\upsilon})]=G_{t,t+\delta}^{t,x;\upsilon}[Y_{t+\delta}^{t,x;\upsilon}].
\end{equation*}
Then we can obtain the DPP for our stochastic optimal control problems as follow:
\begin{theorem}
The value function $u(t,x)$ have the following proposition: for every $0 \leq \delta \leq T-t$, we have
\begin{equation}\label{15}
u(t,x)=\underset{\upsilon(\cdot)\in \mathcal{U}}{sup}\ G_{t,t+\delta}^{t,x;\upsilon}[u(t+\delta,X_{t+\delta}^{t,x;\upsilon})].
\end{equation}
\end{theorem}
\proof
We have
\begin{align*}
u(t,x)=\underset{\upsilon(\cdot)\in \mathcal{U}}{sup}G_{t,T}^{t,x;\upsilon}[\Phi(X_{T}^{t,x;\upsilon})]
=\underset{\upsilon(\cdot)\in \mathcal{U}}{sup}G_{t,t+\delta}^{t,x;\upsilon}[Y_{t+\delta}^{t+\delta,X_{t+\delta}^{t,x;\upsilon};\upsilon}].
\end{align*}
Obviously, $X_{t+\delta}^{t,x;\upsilon}$ is $\mathcal{F}^{B}_{t+\delta}$ measurable. So by Lemma 3.5 and Theorem2.19, we have
\begin{equation*}
u(t,x) \leq \underset{\upsilon(\cdot)\in \mathcal{U}}{sup}G_{t,t+\delta}^{t,x;\upsilon}[u(t+\delta,X_{t+\delta}^{t,x;\upsilon})].
\end{equation*}
Besides, for $\varepsilon >0$, there exists an  admissible control $\bar{\upsilon}(\cdot)\in \mathcal{U}$ such that
\begin{equation*}
u(t+\delta,X_{t+\delta}^{t,x;\upsilon})\leq Y_{t+\delta}^{t+\delta,X_{t+\delta}^{t,x;\upsilon};\bar{\upsilon}}+\varepsilon.
\end{equation*}
Then
\begin{align*}
u(t,x)&\geq \underset{\upsilon(\cdot)\in \mathcal{U}}{sup} G_{t,t+\delta}^{t,x;\upsilon}[u(t+\delta,X_{t+\delta}^{t,x;\upsilon})-\varepsilon]\\
&\geq \underset{\upsilon(\cdot)\in \mathcal{U}}{sup} G_{t,t+\delta}^{t,x;\upsilon}[u(t+\delta,X_{t+\delta}^{t,x;\upsilon})]-C\varepsilon.
\end{align*}
Because $\varepsilon$ can be arbitrarily small, we get \eqref{15}.
\endproof
\begin{proposition}
$u(t,x)$ is $\frac{1}{2}$-$H\ddot{o}lder$ continuous in $t$.
\end{proposition}
\proof
For any given $(t,x)\in [0,T]\times \mathbb{R}^{n}$ and $\delta > 0(t+\delta \leq T)$, from Theorem 3.6, we know that for $\varepsilon > 0$, there exists a $\upsilon(\cdot)\in\mathcal{U}$ such that
\begin{equation*}
G_{t,t+\delta}^{t,x;\upsilon}[u(t+\delta,X_{t+\delta}^{t,x;\upsilon})]+\varepsilon \geq u(t,x) \geq G_{t,t+\delta}^{t,x;\upsilon}[u(t+\delta,X_{t+\delta}^{t,x;\upsilon})].
\end{equation*}
Then we need to prove
\begin{equation}\label{16}
u(t,x)-u(t+\delta,x) \leq C \delta^{\frac{1}{2}}\    (respectively,\ \geq -C \delta^{\frac{1}{2}}).
\end{equation}
We only check the first inequality in \eqref{16}. The second can be proved similarly. We have $\forall \varepsilon > 0$,
\begin{equation}\label{20}
u(t,x)-u(t+\delta,x)\leq I_{\delta}^{1}+ I_{\delta}^{2}+ \varepsilon,
\end{equation}
where
\begin{align*}
I_{\delta}^{1}&=G_{t,t+\delta}^{t,x;\upsilon}[u(t+\delta,X_{t+\delta}^{t,x;\upsilon})]-G_{t,t+\delta}^{t,x;\upsilon}[u(t+\delta,x)], \\ I_{\delta}^{2}&=G_{t,t+\delta}^{t,x;\upsilon}[u(t+\delta,x)]-u(t+\delta,x).
\end{align*}
From Proposition 3.1, we have
\begin{equation*}
\hat{E}_{t}[|X_{t+\delta}^{t.x;\upsilon}-x|^{2}] \leq C(1+|x|^{2}) \delta.
\end{equation*}
By proposition 3.2 and Lemma 3.4, we deduce that
\begin{equation*}
|I_{\delta}^{1}|\leq [C \hat{E}_{t}[|u(t+\delta,X_{t+\delta}^{t,x;\upsilon})-u(t+\delta,x)|^{2}]]^{\frac{1}{2}}
\leq[C\hat{E}_{t}[|X_{t+\delta}^{t,x;\upsilon}-x|^{2}]]^{\frac{1}{2}}\leq C'\delta^{\frac{1}{2}}.
\end{equation*}
Based on the definition of $G_{t,t+\delta}^{t,x;\upsilon}$, we get
\begin{align*}
I_{\delta}^{2}=&\hat{E}_{t}[u(t+\delta,x)+\int_{t}^{t+\delta}f(s,X_{s}^{t,x;\upsilon},Y_{s}^{t,x;\upsilon},
Z_{s}^{t,x;\upsilon},\upsilon_{s})ds\\
&+\underset{i,j=1}{\overset{d}{\sum}}\int^{t+\delta}_{t}g_{ij}(s,X_{s}^{t,x;\upsilon},Y_{s}^{t,x;\upsilon},
Z_{s}^{t,x;\upsilon},\upsilon_{s})d\langle B^{i},B^{j} \rangle_{s}\\
&-\int^{t+\delta}_{t}Z_{s}^{t,x;\upsilon}dB_{s}-(K_{T}^{t,x;\upsilon}-K_{t}^{t,x;\upsilon})]-u(t+\delta,x)\\
=&\hat{E}_{t}[\int_{t}^{t+\delta}f(s,X_{s}^{t,x;\upsilon},Y_{s}^{t,x;\upsilon},Z_{s}^{t,x;\upsilon},\upsilon_{s})ds\\
&+\underset{i,j=1}{\overset{d}{\sum}}\int^{t+\delta}_{t}g_{ij}(s,X_{s}^{t,x;\upsilon},
Y_{s}^{t,x;\upsilon},Z_{s}^{t,x;\upsilon},\upsilon_{s})d\langle B^{i},B^{j} \rangle_{s}]\\
\leq& C' \delta^{\frac{1}{2}}(1+\hat{E}_{t}[\int_{t}^{t+\delta}|X_{s}^{t,x;\upsilon}|^{2}
+|Y_{s}^{t,x;\upsilon}|^{2}+|Z_{s}^{t,x;\upsilon}|^{2}ds]^{\frac{1}{2}}).
\end{align*}
By Proposition 3.2, we can prove the following inequality easily by the similar method in Proposition 3.5 of \cite{r2}
\begin{equation*}
\hat{E}_{t}[\int_{t}^{t+\delta}|Z_{s}^{t,x;\upsilon}|^{2}ds]^{\frac{1}{2}} \leq  C(1+|x|).
\end{equation*}
So we have $I_{\delta}^{2} \leq C' \delta^{\frac{1}{2}}$.
Hence, by \eqref{20} we have
\begin{equation*}
u(t,x)-u(t+\delta,x) \leq C'\delta^{\frac{1}{2}}+\varepsilon.
\end{equation*}
Let $\varepsilon \rightarrow 0$, we  obtain the first inequality of \eqref{16}.
The proof is completed.
\endproof
\section{Value Function and Viscosity Solution of Fully Nonlinear Second-Order Partial Differential Equation}
In this section, we consider the following fully nonlinear second-order partial differential equation
\begin{equation}\label{17}
\left \{
\begin{array}{rl}
&\partial_{t}u+F(D_{x}^{2}u,D_{x}u,u,x,t)=0,\ (t,x)\in[0,T]\times\mathbb{R}^{n},\\
&u(T,x)=\Phi(x),\\
\end{array}
\right.
\end{equation}
where
\begin{align*}
F(D_{x}^{2}u,D_{x}u,u,x,t)=&\underset{\upsilon\in V}{sup}\{G(H(D_{x}^{2}u,D_{x}u,u,x,t,\upsilon))+\langle b(t,x,\upsilon),D_{x}u \rangle\\
&+f(t,x,u,\langle \sigma_{1}(t,x),D_{x}u \rangle,...,\langle \sigma_{d}(t,x),D_{x}u \rangle,\upsilon)\},
\end{align*}
\begin{align*}
H_{ij}(D_{x}^{2}u,D_{x}u,u,x,t,\upsilon)=&\langle D_{x}^{2}u\cdot \sigma_{i}(t,x,\upsilon),\sigma_{j}(t,x,\upsilon) \rangle+2\langle D_{x}u,h_{ij}(t,x,\upsilon) \rangle\\
&+2g_{ij}(t,x,u,\langle \sigma_{1}(t,x,\upsilon),D_{x}u \rangle,...,\langle \sigma_{d}(t,x,\upsilon),D_{x}u \rangle, \upsilon).
\end{align*}
\begin{remark}
The definition and uniqueness of viscosity solution of above second-order partial differential equation can be found in Appendix C in Peng\cite{r3}. So we only need to prove that $u(t,x)$ is a viscosity solution of equation \eqref{17}. Besides, from the result of section 3, we can have that $u(t,x)$ is continuous in $[0,T]\times \mathbb{R}^{n}$.
\begin{definition}
A real-valued continuous function $u(t,x)\in C([0,T]\times\mathbb{R}^{n})$, $u(T,x) \leq \Phi(x)$, for any $x\in \mathbb{R}^{n}$, is called a viscosity sub-solution (super-solution) of \eqref{17}, if   for all functions $\varphi \in C^{2,3}([0,T]\times\mathbb{R}^{n})$ satisfy  $\varphi \geq u$ and $\varphi(t,x)=u(t,x)$ at fixed $(t,x)\in [0,T)\times\mathbb{R}^{n}$ , we have
\begin{equation*}
\partial_{t}\varphi(t,x)+F(D^{2}_{x}\varphi(t,x),D_{x}\varphi(t,x),\varphi(t,x),x,t)\geq 0(\leq 0).
\end{equation*}
\end{definition}

\end{remark}
\begin{theorem}
Under the assumptions (H3) and (H4), the value function $u(t,x)$ defined by \eqref{18} is a viscosity solution of equation \eqref{17}.
\end{theorem}
In order to prove the Theorem, we need three Lemma. Firstly, we set
\begin{align*}
&F_{1}(r,x,y,z,\upsilon)\\
&=\langle b(r,x,\upsilon), D_{x}\varphi(r,x) \rangle+\partial_{t}\varphi(t,x)\\
&~~~~+f(r,x,y+\varphi(r,x),z+( \langle \sigma_{1}(t,x,\upsilon),D_{x}\varphi(r,x)\rangle,...,\langle \sigma_{d}(t,x,\upsilon),D_{x}\varphi(r,x) \rangle),\upsilon),
\end{align*}
\begin{align*}
&F_{2}^{ij}(r,x,y,z,\upsilon)\\
&=\langle D_{x}\varphi(r,x),h_{ij}(r,x,\upsilon) \rangle+\frac{1}{2}\langle D_{x}^{2}\varphi(r,x)\sigma_{i}(r,x,\upsilon),\sigma_{j}(r,x,\upsilon) \rangle\\
&~~~~+g_{ij}(r,x,y+\varphi(r,x),z+(\langle \sigma_{1}(t,x,\upsilon),D_{x}\varphi(r,x) \rangle,...,\langle \sigma_{d}(r,x,\upsilon),D_{x}\varphi(r,x) \rangle),\upsilon).
\end{align*}
Then we consider a G-BSDE defined on the interval $[t,t+\delta] (0 < \delta \leq T-t)$:
\begin{align}\label{19}
Y_{s}^{1,\upsilon}=&\int_{s}^{t+\delta}F_{1}(r,X_{r}^{t,x;\upsilon},Y_{r}^{1,\upsilon},Z_{r}^{1,\upsilon},\upsilon_{r})dr
+\int_{s}^{t+\delta}Z_{r}^{1,\upsilon}dB_{r}-(K_{t+\delta}^{1}-K_{s}^{1})\nonumber\\
&-\underset{i,j=1}{\overset{d}{\sum}}\int_{s}^{t+\delta}F_{2}^{ij}(r,X_{r}^{t,x;\upsilon},Y_{r}^{1,\upsilon},
Z_{r}^{1,\upsilon},\upsilon_{r})d \langle B^{i},B^{j} \rangle_{r},
\end{align}
where $\upsilon(\cdot)\in\mathcal{U}$ and $X_{s}^{t,x;\upsilon}$ defined by \eqref{formula3}.
\begin{lemma}
For $s \in [t,t+\delta]$, we have
\begin{equation*}
G_{s,t+\delta}^{t,x;\upsilon}[\varphi(X_{t+\delta}^{t,x;\upsilon},t+\delta)]-\varphi(X_{s}^{t,x;\upsilon},s)
\end{equation*}
is the solution of \eqref{19}.
\end{lemma}
\proof
From the definition of $G_{s,t+\delta}^{t,x;\upsilon}$, we know that $G_{s,t+\delta}^{t,x;\upsilon}[\varphi(X_{t+\delta}^{t,x;\upsilon},t+\delta)]$ is the solution of G-BSDE \eqref{formula6} on $[t,t+\delta]$ with terminal condition $\varphi(X_{t+\delta}^{t,x;\upsilon},t+\delta)$. Applying It\^{o}'s formula to $\varphi(X_{s}^{t,x;\upsilon},s)$, we can obtain the result.
\endproof

Now we construct a simple G-BSDE by replacing the driving process $X_{s}^{t,x;\upsilon}$ by its deterministic initial value $x$ as follow :
\begin{align}\label{21}
Y_{s}^{2,\upsilon}=&\int_{s}^{t+\delta}F_{1}(r,x,Y_{r}^{2,\upsilon},Z_{r}^{2,\upsilon},\upsilon_{r})dr
+\underset{i,j=1}{\overset{d}{\sum}}\int_{s}^{t+\delta}F_{2}^{ij}(r,x,Y_{r}^{2,\upsilon},Z_{r}^{2,\upsilon},\upsilon_{r})d \langle B^{i},B^{j} \rangle_{r}\nonumber\\
&-\int_{s}^{t+\delta}Z_{r}^{2,\upsilon}dB_{r}-(K_{t+\delta}^{2}-K_{s}^{2}).
\end{align}

\begin{lemma}
We have the following estimate,  for $\upsilon(\cdot)\in \mathcal{U}$,
\begin{equation*}
|Y_{t}^{1,\upsilon}-Y_{t}^{2,\upsilon}| \leq C \delta^{\frac{3}{2}}.
\end{equation*}
Where $C$ is independent of the control processes $\upsilon(\cdot)$.
\end{lemma}
\proof
By proposition 3.1, we have the estimate for $p \geq 2$
\begin{equation*}
\hat{E}_{t}[\underset{s\in [t,t+\delta]}{\sup}|X_{s}^{t,x;\upsilon}-x|^{p}]\leq C(1+|x|^{p})\delta^{\frac{p}{2}}.
\end{equation*}
By proposition 2.18, we get for fixed $p > 2$ and $2 < p < \beta$,
\begin{align*}
|Y_{t}^{1,\upsilon}-Y_{t}^{2,\upsilon}|^{2}&\leq \hat{E}[\underset{s\in[t,t+\delta]}{sup}|Y_{t}^{1,\upsilon}-Y_{t}^{2,\upsilon}|^{2}]\\
&\leq C \{\hat{E}[\underset{s\in[t,t+\delta]}{\sup}\hat{E}_{s}[(\int_{t}^{t+\delta}\hat{F}_{r}dr)^{p}]])^{\frac{2}{p}}
+\hat{E}[\underset{s\in[t,t+\delta]}{\sup}\hat{E}_{s}[(\int_{t}^{t+\delta}\hat{F}_{r}dr)^{p}]]\},
\end{align*}
where
\begin{align*}
\hat{F}_{r}=&|F_{1}(r,X_{r}^{t,x;\upsilon},Y_{r}^{2,\upsilon},Z_{r}^{2,\upsilon},\upsilon_{r})
-F_{1}(r,x,Y_{r}^{2,\upsilon},Z_{r}^{2,\upsilon},\upsilon_{r})|\\
&+\underset{i,j=1}{\overset{d}{\sum}}|F^{i,j}_{1}(r,X_{r}^{t,x;\upsilon},
Y_{r}^{2,\upsilon},Z_{r}^{2,\upsilon},\upsilon_{r})-F^{i,j}_{2}(r,x,Y_{r}^{2,\upsilon},Z_{r}^{2,\upsilon},\upsilon_{r})|.
\end{align*}
It is easy to prove that
\begin{equation*}
\hat{F}_{r} \leq C |X_{r}^{t,x;\upsilon}-x|.
\end{equation*}
Then we can deduce that
$|Y_{t}^{1,\upsilon}-Y_{t}^{2,\upsilon}| \leq C \delta^{\frac{3}{2}}$.

\endproof
\begin{lemma}
We have
\begin{equation*}
\underset{\upsilon(\cdot)\in \mathcal{U}}{sup}Y_{t}^{2,\upsilon}=Y^{0}(t).
\end{equation*}
where $Y_{0}(\cdot)$ is the solution of the following ODE:
\begin{equation*}
\left \{
\begin{array}{rl}
&-dY^{0}_{s}=F^{0}(s,x,Y^{0}_{r},0)ds,\ s\in[t,t+\delta],\\
&Y^{0}_{t+\delta}=0,\\
\end{array}
\right.
\end{equation*}
where $F^{0}(r,x,y,z)=\underset{\upsilon\in V}{sup}\{F_{1}(r,x,y,z,\upsilon)+2G[(F_{2}^{ij}(r,x,y,z,\upsilon))_{i,j=1}^{d}]\}$.
\end{lemma}
\proof
By Theorem 2.16, we know that the G-BSDE \eqref{21} have a unique solution $(Y,Z,K)$. Hence there exists a process
\begin{equation*}
 V^{2,\upsilon}_{s}=\underset{i,j=1}{\overset{d}{\Sigma}}\int_{t}^{s}F_{2}^{ij}(r,x,Y^{2,\upsilon}_{r},Z^{2,\upsilon}_{r},\upsilon)d \langle B^{i},B^{j} \rangle_{r}-\int_{t}^{s}2G((F_{2}^{ij}(r,x,Y^{2,\upsilon}_{r},Z^{2,\upsilon}_{r},\upsilon))_{i,j=1}^{d})dr.
\end{equation*}
Here $V_{s}^{2,\upsilon}$, $s\in[t,t+\delta]$ is a decreasing and continuous process by \cite{r38}. Besides, it satisfies   $\hat{E}[\underset{s\in[t,t+\delta]}{\sup}|V_{s}^{2,\upsilon}|^{\beta}] < \infty$ obviously.
So $Y^{2,\upsilon}_{s}$ is the solution of the following G-BSDE:
\begin{align*}
Y^{2,\upsilon}_{s}=&\int_{s}^{t+\delta}[F_{1}(r,x,Y^{2,\upsilon}_{r},Z^{2,
\upsilon}_{r},\upsilon_{r})+2G[(F_{2}^{ij}(r,x,Y^{2,\upsilon}_{r},Z^{2,\upsilon}_{r},\upsilon_{r}))_{i,j=1}^{d}]dr\\
&-\int_{s}^{t+\delta}Z_{r}^{2,\upsilon}dB_{r}
-(K_{t+\delta}^{2}-K_{s}^{2})+V^{2,\upsilon}_{t+\delta}-V^{2,\upsilon}_{s},
\end{align*}
where $\upsilon(\cdot)\in \mathcal{U}$.
In addition, we have
\begin{equation*}
Y^{0}_{t}=\int_{s}^{t+\delta}F^{0}(r,x,Y^{0}_{r},Z^{0}_{r})dr- \int_{s}^{t+\delta}Z^{0}_{r}dB_{r}-(K^{0}_{t+\delta}-K^{0}_{s})+(V^{0}_{t+\delta}-V^{0}_{s}),
\end{equation*}
where $(Z,K,V)=0$.
By the comparison theorem 2.20 and the definition of $F^{0}$, we have for $\upsilon(\cdot)\in \mathcal{U}$,
\begin{equation*}
Y_{s}^{2,\upsilon} \leq Y_{s}^{0}, s\in[t,t+\delta].
\end{equation*}
On the other hand, there exists a measurable function $\upsilon'(r,x,y,z):[t,T]\times\mathbb{R}^{n}\times\mathbb{R}\times\mathbb{R}^{d}\times\mathbb{R}\rightarrow V$ such that
\begin{equation*}
F^{0}(r,x,y,z)=F_{1}(r,x,y,z,\upsilon')+2G[(F_{2}^{ij}(r,x,y,0,\upsilon'))_{i,j=1}^{d}].
\end{equation*}
Then we have $\upsilon'(r,x,Y_{r}^{0},Z_{r}^{0}) \in \mathcal{U}$ and $Y^{0}_{t}$ is the solution of following G-BSDE:
\begin{align*}
Y^{0}_{s}=&\int_{s}^{t+\delta}F_{1}(r,x,Y_{r}^{0},Z_{r}^{0},
\upsilon'_{r})dr+\underset{i,j=1}{\overset{d}{\sum}}\int_{s}^{t+\delta}F_{2}^{ij}(r,x,Y_{r}^{0},Z_{r}^{0},\upsilon'_{r})d \langle B^{i},B^{j} \rangle_{r}\\
&- \int_{s}^{t+\delta}Z^{0}_{r}dB_{r}-(K^{0}_{t+\delta}-K^{0}_{s}),
\end{align*}
where $Z^{0}_{r,\upsilon}=0$,
\begin{equation*}
K_{s}^{0}=\underset{i,j=1}{\overset{d}{\sum}}\int_{t}^{s}F_{2}^{ij}(r,x,Y^{0}_{r},0,\upsilon')d \langle B^{i},B^{j} \rangle_{r}-\int_{t}^{s}2G((F_{2}^{ij}(r,x,Y^{0}_{r},0,\upsilon'))_{i,j=1}^{d})dr.
\end{equation*}
So $Y^{0}_{t}\leq \underset{\upsilon(\cdot)\in \mathcal{U}}{sup}Y_{t}^{2,\upsilon}$. Now we have proved the lemma.

\endproof
Then we give the proof of Theorem 4.3:
\proof
We set $\varphi\in C^{2,3}([0,T]\times \mathbb{R}^{n})$ and $\varphi(t,x)=u(t,x)$ for fixed $(t,x)\in[0,T]\times \mathbb{R}^{n}$. From Theorem 3.6, we know
\begin{equation*}
\varphi(t,x)=u(t,x)=\underset{\upsilon(\cdot)\in\mathcal{U}}{sup}G_{t,t+\delta}^{t,x;\upsilon}[u(X_{t+\delta}^{t,x},t+\delta)].
\end{equation*}
By $\varphi \geq u(\varphi \leq u)$ and the definition of $G$
\begin{equation*}
\underset{\upsilon(\cdot)\in\mathcal{U}}{sup}\{G_{t,t+\delta}^{t,x;\upsilon}[u(X_{t+\delta}^{t,x},t+\delta)]-\varphi(t,x)\}\geq 0(\leq 0).
\end{equation*}
Then form lemma 4.4
\begin{equation*}
\underset{\upsilon(\cdot)\in\mathcal{U}}{sup}Y_{t}^{1,\upsilon}\geq 0(\leq 0).
\end{equation*}
Besides from lemma 4.5
\begin{equation*}
\underset{\upsilon(\cdot)\in\mathcal{U}}{sup}Y_{t}^{2,\upsilon}\geq C\delta^{\frac{3}{2}}(\leq C\delta^{\frac{3}{2}}).
\end{equation*}
Finally, lemma 4.6 implies
\begin{equation*}
Y^{0}(t)\geq C\delta^{\frac{3}{2}}(\leq C\delta^{\frac{3}{2}}).
\end{equation*}
So $F^{0}(r,x,0,0)\geq 0(\leq 0)$ and from the definition of viscosity solution of equation \eqref{17}, we know $u(t,x)$ is a viscosity solution of equation \eqref{17}.
\endproof

\end{document}